# Importance of Linking Inertia and Frequency Response Procurement: The Great Britain Case

Aimon Mirza Baig, Luis Badesa and Goran Strbac
Department of Electrical and Electronic Engineering
Imperial College London
London SW7 2AZ, United Kingdom
Email: {a.mirza-baig19, luis.badesa, g.strbac}@imperial.ac.uk

*Abstract*—In order to decarbonise the electricity sector, the future Great Britain (GB) power system will be largely dominated by non-synchronous renewables. This will cause low levels of inertia, a key parameter that could lead to frequency deterioration. Therefore, the requirement for ancillary services that contain frequency deviations will increase significantly, particularly given the increase in size of the largest possible loss with the commissioning of large nuclear plants in the near future. In this paper, an inertia-dependent Stochastic Unit Commitment (SUC) model is used to illustrate the benefits of linking inertia and frequency response provision in low-inertia systems. We demonstrate that the cost of procuring ancillary services in GB could increase by 165% if the level of inertia is not explicitly considered when procuring frequency response. These results highlight the need to re-think the structure of ancillary-services markets, which in GB are nowadays held one month ahead of delivery.

*Keywords*—*Ancillary services, frequency stability, renewable energy, unit commitment.*

## Nomenclature

*Indices and sets*

| | |
|---|---|
| $g$, $G$ | Index and set of thermal generators |
| $n$, $N$ | Index and set of nodes in scenario tree |
| $s$, $S$ | Index and set of storage units |

*Constants*

| | |
|---|---|
| $\Delta f_{max}$ | Maximum admissible frequency deviation at the nadir (Hz) |
| $\Delta \tau(n)$ | Time interval corresponding to node $n$ (h) |
| $\pi(n)$ | Probability of reaching node $n$ in the SUC |
| $c_g^m$ | Marginal cost of thermal unit $g$ (£/MWh) |
| $c_g^{nl}$ | No-load cost of thermal unit $g$ (£/h) |
| $c_g^{st}$ | Start-up cost of thermal unit $g$ (£) |
| $f_0$ | Nominal frequency of the power grid (Hz) |
| $H_g$ | Inertia constant of thermal unit $g$ (s) |
| $P_g^{max}$ | Maximum generation of thermal unit $g$ (MW) |
| $RoCoF_{max}$ | Maximum Rate-of-Change-of-Frequency admissible (Hz/s) |
| $T_{PFR}$ | Ramping time for PFR (s) |
| $T_{EFR}$ | Ramping time for EFR (s) |

*Decision variables*

| | |
|---|---|
| $P_L$ | Largest power infeed (MW) |
| $R_g$ | PFR provision from thermal unit $g$ (MW) |
| $R_s$ | EFR provision from storage unit $s$ (MW) |
| $P_g(n)$ | Power output of thermal unit $g$ at node $n$ (MW) |

*Linear expressions*

| | |
|---|---|
| $C_g(n)$ | Operating cost of thermal unit $g$ at node $n$ in the SUC (£) |
| EFR | Enhanced Frequency Response (MW) |
| $H$ | System inertia (GW·s) |
| $N_g^{sg}(n)$ | Number of thermal units $g$ that start generating at node $n$ in the SUC |
| $N_g^{up}(n)$ | Number of thermal units $g$ that are online at node $n$ in the SUC |
| PFR | Primary Frequency Response (MW) |

## I. Introduction

Increasing awareness by the public on the adverse effects of greenhouse gasses on the environment, has led institutions and regulators to promote the use of renewable energy sources for electricity generation. The British parliament passed a legally binding commitment to achieve net-zero emissions by 2050. To achieve the above target, it is necessary to decarbonise the UK power sector, as this sector contributes to roughly 25% of total emissions in the UK [1].

In the past, the UK power sector was dominated by synchronous thermal generators such as coal, nuclear and gas power plants. However, the future UK generation mix will be dominated by renewables (predominantly wind and solar PV) combined with a nuclear fleet composed by a small number of large power plants. Decarbonisation of the power sector involves considerable challenges due to the uncertain nature of renewables and low levels of system inertia induced by non-synchronous generators, which causes frequency deterioration [2, 3].

The low level of inertia in the system significantly impacts the frequency deterioration that follows an energy-demand imbalance caused by the sudden loss of a generator. To contain this frequency deterioration, provision of auxiliary services in the system is required. Therefore, understanding the value of inertia is necessary to inform market operators on best practices for market design for ancillary services. This paper will focus on understanding the importance of linking inertia and response procurement in low-inertia systems, particularly in electricity grids with a large size of the worst possible contingency.

Research has been carried out on analysing the economic value of inertia. The authors in [4] quantified the benefits of increased system inertia to limit the Rate-of-Change-of-Frequency (RoCoF) after a generation loss, while the work in [5] developed a marginal-pricing scheme that assigns a price to inertia depending on the system condition. Reference [6] analyses the value of inertia and demonstrates that this value varies between hours. This paper used a generation scheduling model that respects frequency stability to quantify the value of inertia, and highlighted that this value depends on the system operating condition at any given hour. However, the aforementioned works do not analyse the value of linking inertia to the frequency response market. Therefore, the present paper focuses on understanding the need for linking inertia and frequency response procurement, as these two services both provide frequency support although with different dynamics. The current approach for procuring ancillary services for frequency support in Great Britain (GB) is replicated here, where the system operator National Grid unlinks the provision of inertia and frequency response. Frequency response is procured in month-ahead auctions, while the level of inertia can only be known with acceptable accuracy in the day-ahead energy market [7]. This approach is compared to a full co-optimisation of inertia and response, where the energy and ancillary services markets are cleared simultaneously and therefore the optimal balance between these two services can be found.

The inertia-dependent, frequency-security constrained Stochastic Unit Commitment (SUC) model developed in [8] has been used for this work, which assures frequency security by optimally scheduling inertia and frequency response. This model also allows to account for the uncertainty in a system with high penetration of Renewable Energy Sources (RES), by scheduling energy production and reserves while minimising the operational cost. The key contributions of this work are:

1) The importance of linking inertia and frequency response procurement is demonstrated through several relevant case studies.
2) The results presented can inform the design of ancillary-services markets that consider the interaction between inertia and response in order to reduce the overall system cost.
3) Both the current and future GB system are analysed, highlighting that linking these services will become increasingly important with the increment in the largest loss to 1.8GW after the commissioning of nuclear plant Hinkley Point C [9] and increased wind penetration to meet emissions targets.

The rest of this paper is organised as follows. Section II describes the unit commitment model with frequency constraints used to conduct simulations, whereas section III explains the dynamics of the post-fault frequency requirements. Section IV explains the results obtained in terms of economic benefits of co-optimising inertia and frequency response. Finally, section V provides the conclusion.

## II. UNIT COMMITMENT WITH FREQUENCY CONSTRAINTS

The optimal selection of a combination of thermal units to be committed in order to meet the electric demand is called the Unit Commitment (UC) problem. Traditionally, a deterministic UC has been used as the tool used to schedule power system's generation, given that uncertainty was limited to small variations around the demand forecast. In addition to energy, the UC can also be used to schedule reserves required to compensate a possible mismatch between generation and demand due to deviations from the demand and RES forecasts.

However, the increasing RES penetration in most power systems has reduced the effectiveness of the deterministic techniques used for scheduling. Deterministic models can either increase the operational cost of the system due to over-scheduling of reserves, or increase the risk of load shed. Research in [6] shows that the stochastic counterpart of UC can reduce the operational cost of a system by up to 4% compared to the deterministic UC approach, for systems with high RES penetration. This reduction in the cost is due to optimally scheduling the combination of spinning reserves from Combined Cycle Gas Turbines (CCGTs) and standing reserves from Open Cycle Gas Turbines (OCGTs). Although CCGTs have lower operational cost, they are less flexible due to longer start-up times, whereas OCGTs have higher operational cost but are highly flexible with almost zero start-up time.

In this work, we use the SUC with frequency-stability constraints developed in [8] which also schedules inertia and frequency response. This SUC is used to quantify the economic value of linking inertia and response, as this model considers uncertainty in RES generation and seeks to optimally schedule energy production and delivery of ancillary services whilst minimising the expected operation cost. The infinite set of possible realisations of RES generation are modelled in the form of a scenario tree, which is built by a quantile-based scenario method as in [10].

A rolling planning approach is considered in the SUC in this paper, which is based on two steps: first, the whole optimisation is carried out for a 24-hour period in an hourly time step and the decision in the current-time node is considered and applied, while the rest of the decisions are discarded. Then, in the next time step the realisation of the stochastic variable and updated forecasts become available and then a new scenario tree is built for the next 24 hours. This building of scenario trees for each time step is repeated for the entire duration of the simulations, therefore this methodology is called 'rolling planning'. The objective function of the SUC is to minimise the operational cost of the system, and it is given by:

$$\min \ \sum_{n \in N} \pi(n) \sum_{g \in G} C_g(n) \qquad (1)$$

Where the operational cost of the thermal generators is given by the following equation:

$$C_g(n) = c_g^{st} N_g^{sg}(n) + \Delta\tau(n)[c_g^{nl} N_g^{up}(n) + c_g^m P_g(n)] \qquad (2)$$

The objective function of the SUC model is subject to the frequency-stability constraints, which are discussed in next section of this paper.

## III. REQUIREMENTS FOR POST-FAULT FREQUENCY DYNAMICS

To guarantee a secure supply of electricity, frequency in the electricity grid must be within a narrow band around the nominal value (50Hz in GB). Frequency deviation outside the secure levels could damage the different devices in the grid, and therefore automatic generation and/or demand disconnection could be triggered by protection devices. The electric frequency in the grid is determined by the rotating speed of the synchronous generators, which spontaneously slow down when demand is higher than generation (for example, after a generation outage) and vice versa. This slowing down causes a frequency deviation, which must be contained to guarantee the stability in the grid.

In order to maintain the frequency within the secure range, the system operator must schedule enough frequency response and inertia to cover for the loss of any generator in the system. Frequency response has traditionally been provided in GB by the Primary Frequency Response (PFR) service, which is required to ramp up by 10s after a fault and is typically provided by thermal generators. However, to tackle the low-inertia problem, National Grid introduced in 2017 a new service called Enhanced Frequency Response (EFR) [11], which is required to ramp up in just 1s and is currently provided by battery storage units.

To calculate the volume of inertia and response that must be scheduled at any given time, the conditions that guarantee frequency to stay within secure limits can be deduced from the swing equation, which represents the dynamics of frequency after a fault [12].

$$\frac{2H}{f_o} \cdot \frac{d\Delta f(t)}{dt} = \sum_{s \in S} \text{EFR}_s(t) + \sum_{g \in G} \text{PFR}_g(t) - P_L \quad (3)$$

where EFR and PFR are modelled as:

$$\text{EFR}_s(t) = \begin{cases} \frac{R_s}{T_{\text{EFR}}} \cdot t & \text{if } t \leq T_{\text{EFR}} \\ R_s & \text{if } t > T_{\text{EFR}} \end{cases} \quad (4)$$

$$\text{PFR}_g(t) = \begin{cases} \frac{R_g}{T_{\text{PFR}}} \cdot t & \text{if } t \leq T_{\text{PFR}} \\ R_g & \text{if } t > T_{\text{PFR}} \end{cases} \quad (5)$$

The three limits that must be respected in any grid are the maximum admissible RoCoF, minimum acceptable value for a frequency drop (called frequency nadir) and the quasi-steady-state (q-s-s) requirement. The RoCoF limit can be obtained from equation (3) and is given by the following inequality:

$$H \geq \frac{P_L \cdot f_o}{2 \, \text{RoCoF}_{max}} \quad (6)$$

The frequency nadir must be above a certain level to avoid the activation of load shedding, and is given by the following expression [8]:

$$\left(\frac{H}{f_o} - \frac{\sum_{s \in S} \text{EFR}_s \cdot T_{\text{EFR}}}{4 \Delta f_{max}}\right) \cdot \sum_{g \in G} \text{PFR}_g \\ \geq \frac{(P_L - \sum_{s \in S} \text{EFR}_s)^2 \cdot T_{\text{PFR}}}{4 \Delta f_{max}} \quad (7)$$

The q-s-s requirement defines a stable value that frequency must reach in 60 seconds after the contingency occurs. The constraint for q-s-s can be directly obtained from the swing equation by considering RoCoF to be effectively zero:

$$\sum_{s \in S} \text{EFR}_s + \sum_{g \in G} \text{PFR}_g \geq P_L \quad (8)$$

To meet the above frequency requirements, inertia, EFR and PFR are considered as decision variables in the SUC, and a sufficient amount of these services is scheduled by the model so that the frequency constraints are respected.

## IV. ECONOMIC BENEFITS OF CO-OPTIMISING INERTIA AND FREQUENCY RESPONSE

This section analyses the importance of considering the system level of inertia when clearing the market for frequency response in low-inertia systems. To do so, we quantify the increase in overall system cost if these two services are unlinked. Given that both these services help contain a frequency drop after a generation-demand imbalance, an optimal balance in volumes of inertia and response scheduled can be found by using the constraints presented in Section III. The difference between inertia and response lies in their different dynamics (i.e. inertia support is proportional to the derivative of frequency while response support is proportional to frequency deviation), and this difference is explicitly considered in the frequency constraints.

To conduct this analysis, the GB system is used as the simulation platform, given that this system is already experiencing low levels of inertia and therefore it represents a particularly illustrative case. The characteristics of thermal plants considered in this work are included in Table I, while the wind capacity and largest loss are the two key sensitivities presented in the case studies. A 2.6GW pump-hydro storage with 10GWh duration and 75% round efficiency is also modelled in the system, as well as 250MW of battery storage with 1GWh duration and 96% round efficiency. The battery storage is considered to have a 200MW capacity for providing EFR. According to GB grid standards, the frequency limits were set to $\text{RoCoF}_{max}$=0.5Hz/s and $\Delta f_{max}$=0.8Hz.

To understand the importance of jointly procuring inertia and response, the current approach in place in GB (which unlinks the inertia and frequency response procurement) is compared to a strategy in which these two services are explicitly co-optimised. These two approaches were carried out under different sensitivities of wind capacity and sizes of the largest loss:

- Unlinked approach for current system ('Unlink 1'): this case replicates the current method used by National Grid, consisting on procuring frequency response through auctions held one month ahead of actual delivery, therefore impeding to accurately estimate the level of inertia when setting the volume of response that will be cleared in the auction. In addition, the current situation of the GB system is considered in this case, by using a 25GW wind capacity and a value of the largest possible generation loss of 1.32GW (driven by the largest nuclear plant currently in the country, named Sizewell B) [13].

- Unlinked approach for future system ('Unlink 2'): this case uses the same approach for procurement of inertia and response as 'Unlink 1', but considering a largest loss of 1.8GW (corresponding of the power rating of nuclear station Hinkley Point C, currently under construction and expected to be commissioned in coming years) and a wind capacity of 50GW, as the share of renewables increases to meet emission targets.

- Co-optimised approach for current system ('Co-opt 1'): this case represents the current GB situation with largest loss of 1.32GW and 25GW wind capacity, but the scheduling of inertia and frequency response is fully co-optimised by jointly scheduling energy and frequency services in the SUC with frequency constraints, as presented in Sections II and III.
- Co-optimised approach for future system ('Co-opt 2'): this case is equivalent to 'Co-opt 1', while considering a largest loss of 1.8GW and wind capacity of 50GW.

TABLE I
CHARACTERISTICS OF THERMAL PLANTS

|  | Nuclear | CCGT | OCGT |
|---|---|---|---|
| Number of Units | 4 | 100 | 30 |
| Rated Power (MW) | 1800 | 500 | 100 |
| Min Stable Generation (MW) | 1800 | 250 | 50 |
| No-Load Cost (£/h) | 0 | 7809 | 8000 |
| Marginal Cost (£/MWh) | 10 | 47 | 200 |
| Start-up Cost (£) | n/a | 10000 | 0 |
| Start-up Time (h) | n/a | 4 | 0 |
| Min Up Time (h) | n/a | 4 | 0 |
| Min Down Time (h) | n/a | 1 | 0 |
| Inertia Constant (s) | 5 | 5 | 5 |
| Max Response (MW) | 0 | 50 | 20 |
| Response Slope | 0 | 0.5 | 0.5 |

The results for these simulations are presented in following sections, highlighting the key aspects that should be paid attention to when designing markets for frequency services.

*A. Importance of co-optimisation in the presence of EFR*

Here we demonstrate how the unlinking of inertia and response can significantly increase system costs in a high renewable system, even swiping away the benefits of 200MW of EFR (a highly valuable service in low-inertia grids given its fast speed). A simulation for the month of January, with an average hourly demand of 43GWh, is considered here. The cost for unlinking and co-optimising inertia and response under different scenarios is computed, considering the following strategies for procurement of frequency services: 1) co-optimised inertia and PFR, for a system with no EFR present; 2) unlinked inertia and PFR, with a fixed volume of 200MW of EFR, as per current practice in GB by 2020 [14]; and 3) co-optimised inertia, PFR and EFR.

The results from the simulations are presented in Fig. 1. By comparing the first and second set of columns, it is clear that the system cost for the 'unlinked' cases is higher, even though a 200MW volume of the valuable EFR is present. The increase in cost is due to the conservativeness introduced by unlinking inertia and response: higher volumes of response are needed in this unlinked case to guarantee system security, leading to higher system costs. This effect is further investigated in the next section.

Fig. 1 also shows that the overall system cost decreases in the future system for all the three cases. This is due to the decrease in total fuel costs for energy production, given the higher wind capacity that can cover a higher share of demand.

Finally, by comparing the first and third set of columns we can quantify the savings that can be drawn from co-optimising not only inertia and PFR, but also 200MW of EFR. The savings are of £10m/month for the current system and of £60m/month for the future system, illustrating the value of the fast service EFR, particularly as the system level of inertia decreases in the future with the increase in non-synchronous wind capacity.

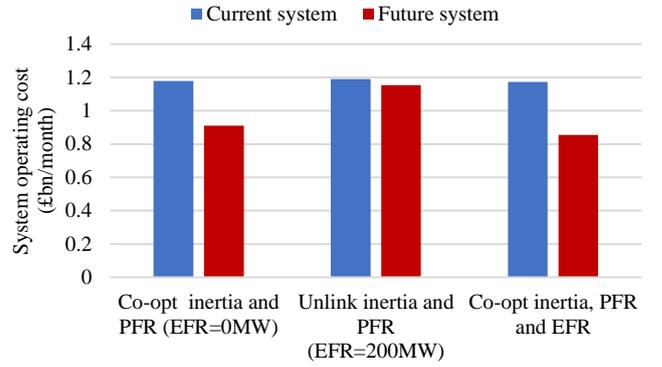

Fig 1. Comparison of monthly system cost under two system scenarios and different strategies for the procurement of frequency services.

The results presented in this section demonstrate the importance of co-optimising inertia and response in order to take full advantage of the lower energy costs from RES, as well as to take full advantage of the fast speed of new response services such as EFR.

*B. Current GB system*

This section particularly focuses on analysing the current approach for procuring response in GB, where the system operator holds monthly auctions for which response is procured up to one month ahead of delivery. To replicate this approach a two-step process is used in this paper: first, an energy-only SUC simulation is run (i.e. an SUC without the frequency constraints implemented) to understand the volume of inertia that results as a by-product of energy production. Then, since inertia is currently not explicitly procured as a distinct service, we take the floor of inertia in the energy-only simulation for each month, and compute the volume of response necessary to respect the frequency constraints in Section III. This computed value for the volume of response needed is then enforced in the SUC for the whole month, and the two-step process is repeated for each month of the year.

The unlinked strategy is compared to a fully co-optimised strategy for each month, considering the characteristics of the current GB power system. The results are presented in Fig. 2, which shows that the co-optimised approach provides lower system cost than the unlinked approach for every month of the year. However, the savings from co-optimisation are not very significant. This result shows why the linking of inertia and response has not been a concern in the past: the size of the largest loss is not very high (1.32GW, compared to the expected 1.8GW in coming years), which results in a moderate requirement for frequency services as defined by the constraints presented in Section III. Moreover, moderate penetration of wind in the system implies that a higher number of synchronous plants are typically online to produce energy, which results in a relatively high level of inertia that reduces the need for response. Nevertheless, the results in the next section of this paper will demonstrate the great importance of co-optimising the different frequency services in the future GB power system.

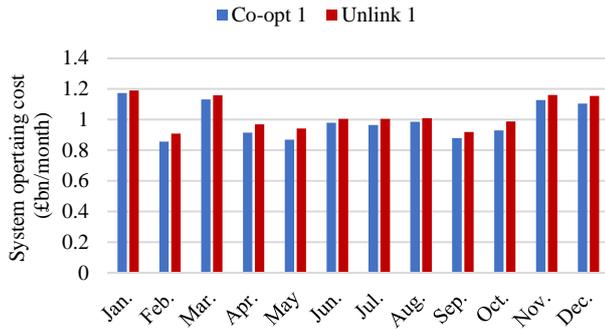

Fig 2. Comparison of monthly system cost from the co-optimised and unlinked strategies for procurement of frequency services, while considering the current GB system.

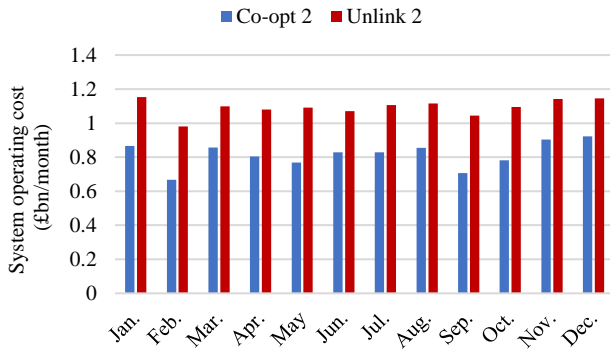

Fig 3. Comparison of monthly system cost from the co-optimised and unlinked strategies for procurement of frequency services, while considering the future GB system.

*C. Future GB system*

Here we consider the future GB system with increased wind capacity and size of the largest loss after the commissioning of plant Hinkley Point C. The two operational strategies, 'co-optimise' and 'unlink', were analysed.

The results in Fig. 3 show a significant increase in the system cost for every month when comparing the 'co-optimise' approach with the 'unlink' approach. This increase is driven by both the high value of the largest loss, which increases the requirement for ancillary services (both inertia and response), and the high wind penetration, which displaces thermal plants to provide energy, thus reducing the level of inertia (which must be compensated with a higher volume of frequency response). Therefore, in a system that exhibits high RES capacity and a large loss to secure against, it is of great importance to link inertia and response procurement in order to achieve a cost-effective operation.

*D. Annual system cost under the different scenarios*

This section presents the annual system cost under the four different scenarios described in the introduction of Section IV. The graph in Fig. 4 shows the annual system cost under these scenarios, which summarises the results presented in Sections IV-B and IV-C. The results again highlight that the co-optimised approach is more beneficial for the future system with increased largest loss and wind capacity (achieving savings of £3.4bn/year), although this strategy is certainly still beneficial for the current system (savings of £0.5bn/year).

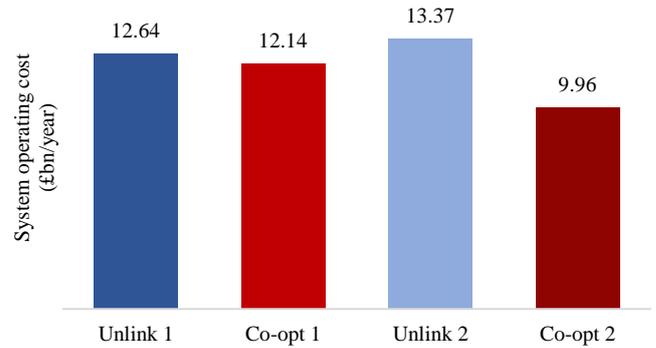

Fig 4. Comparison of annual system cost under the four scenarios.

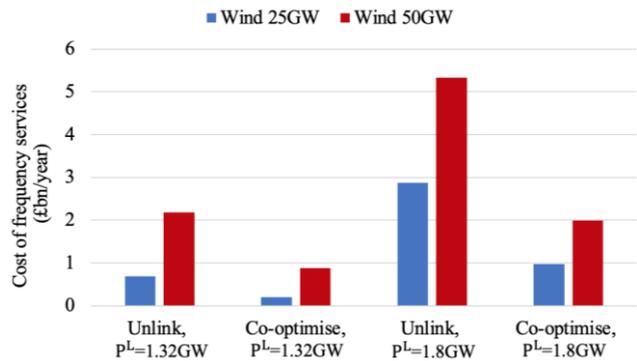

Fig 5. Sensitivity analysis for the cost of frequency services in a system with varying wind capacity and size of the largest possible loss.

As explained in Section IV-C, the higher savings from co-optimisation in the future system are due to both the higher size of the largest loss, which requires a higher amount of inertia and response to keep the system stable, and the higher wind capacity, which has the effect of decreasing the inertia that results as a by-product of energy. In order to analyse which of these two factors is the main cause making the co-optimisation strategy very beneficial, in the next section we further study each of these factors in isolation.

*E. Investigating the main driver for frequency services*

To understand which is the main driver of a higher need for frequency services (a higher loss, or a higher penetration of non-synchronous RES), here we consider the following two 'synthetic' cases (that is, variations on the future GB system that are not expected to occur): 1) the largest loss increases from 1.32GW to 1.8GW, while the wind capacity stays at 25GW; and 2) the wind capacity increases from 25GW to 50GW but the largest loss stays at 1.32GW. Both the co-optimised and unlinked strategies for procurement of inertia and response were considered for each of these cases.

Fig. 5 presents the results of this sensitivity analysis. The cost of frequency services shown in this figure has been obtained by running two SUC simulations, one with the frequency constraints introduced in Section III implemented, and the other without these constraints. The cost of frequency services is computed by taking the difference in the cost obtained in the first simulation minus the cost obtained in the second simulation, as this difference corresponds purely to the cost of frequency services (since the first simulation corresponds to an energy-only SUC).

The results in Fig. 5 show that both the size of the largest loss and the wind capacity have a similar impact on the value of the 'co-optimised' vs. 'unlink' strategies, although the largest loss has a slightly higher impact: the savings from the 'co-optimise' strategy compared to the 'unlink', for a loss of size 1.8GW are of £3.3bn/year with wind capacity of 50GW, but still of £1.9bn/year with a wind capacity of 25GW. On the other hand, the savings for a case of 50GW wind and 1.32GW loss are of £1.2bn/year, lower than the £1.9bn/year mentioned before.

Furthermore, the results also demonstrate that the cost for frequency services can increase by 165% in the expected future system (1.8GW loss, 50GW wind) if inertia and response are unlinked. To further understand why the co-optimised strategy becomes so effective in this system, let's consider the volumes of inertia and response that would be procured in the 'unlinked' case: the amount of PFR required was calculated by taking the floor of inertia for each month from the energy-only SUC, as explained in Section IV-B; since the floor of inertia showed to be very low for every month in this high-wind system, the floor of inertia that needs to be enforced was then determined by the RoCoF constraint eq. (6), which sets a minimum admissible inertia of 90GW·s. This value of inertia then requires a volume of 4.6GW of PFR to meet the nadir constraint eq. (7). In order to procure 4.6GW of PFR in the system, which is provided by the CCGTs that are committed, a high number of these thermal plants is needed, since each of them only has a 50MW capacity of PFR, as described in Table I. In particular, 92 CCGTs would be needed to cover this PFR requirement, which in turn results in an inertia of 230GW·s (since the inertia constant of every plant is assumed to be of 5s). In conclusion, the 'unlink' case leads to an over-procurement of inertia of a factor 2.5, which causes an unnecessary increase in system costs that is particularly relevant in the future GB system.

## V. Conclusion

This paper focuses on analysing the value of co-optimising the procurement of inertia and frequency response services. It has been demonstrated that the cost of procuring these ancillary services could increase by up to 165% in the future GB power system, if the current approach of unlinking the procurement of inertia and response continues to be used. This increase is impacted by both the higher expected value of the largest possible loss (with the commissioning of a large nuclear plant in the near future) and the expected higher capacity of non-synchronous RES to meet emissions targets.

Currently, National Grid is trialling a day-ahead market for frequency response, with the goal of opening the market to new participants such as renewable generators. While increasing the number of participants in the market will certainly have benefits in terms of increased liquidity, the results presented in this paper demonstrate that it could have an additional key benefit: it would allow to link inertia and response procurement. The case studies conducted in this paper have demonstrated that the current approach, consisting on procuring response through month-ahead auctions, would significantly increase the cost of ancillary services in the future power system with high penetration of wind and increased size of the largest loss. Hence, co-optimising inertia and response has proven to be key for a cost-effective operation of the future low-carbon system.

ACKNOWLEDGMENT

This work has been partly supported by EDF Energy R&D UK through project 'Integrated Development of Low-Carbon Energy Systems' (IDLES, grant EP/R045518/1).